\documentclass{amsart}
\usepackage{amssymb}
\usepackage[english,francais]{babel}

\newtheorem{theorem}{Theorem}[section]
\newtheorem{lemma}[theorem]{Lemma}
\newtheorem{e-proposition}[theorem]{Proposition}
\newtheorem{corollary}[theorem]{Corollary}
\newtheorem{e-definition}[theorem]{Definition\rm}

\newtheorem{theoreme}{Th\'eor\`eme}[section]

\newtheorem{proposition}[theoreme]{Proposition}

\def\og{\leavevmode\raise.3ex\hbox{$\scriptscriptstyle\langle\!\langle$~}}
\def\fg{\leavevmode\raise.3ex\hbox{~$\!\scriptscriptstyle\,\rangle\!\rangle$}}

\begin{document}
\title{$S$-adic conjecture and Bratteli diagrams}




\author{Fabien Durand},
\email{fabien.durand@u-picardie.fr}
\author{Julien Leroy}
\email{julien.leroy@u-picardie.fr}

\address{
Universit\'e de Picardie Jules Verne,
Laboratoire Ami\'enois de Math\'ematiques Fondamentales et Appliqu\'ees,
CNRS-UMR 7352,
33 rue Saint Leu,
80039 Amiens Cedex,
France}



\maketitle

\begin{abstract}
In this note we apply a substantial improvement of a result of S. Ferenczi on $S$-adic subshifts to give Bratteli-Vershik representations of these subshifts.
\end{abstract}

\section{Introduction}

In their seminal papers \cite{MH1,MH2} G. A. Hedlund and M. Morse proved that a sequence $x\in A^\mathbb{N}$ is ultimately periodic if and only if $p_x (n) = n$ for some $n$ where $p_x (n)$ is the word complexity of $x$, that is the number of distinct words of length $n$ in $x$. 
Moreover they showed that sequences satisfying $p_x(n)=n+1$ for all $n$ exist, are uniformly recurrent, intimately related to the rotations on the torus and $S$-adic, i.e., produced by an infinite product of finitely many morphisms  ($3$ in fact).
Then, certainly induced by the (sub-affine complexity) examples in \cite{AR} and the complexity of substitutions (\cite{Qu}),
B. Host conjectured that there exists a strong notion of $S$-adicity which is equivalent to sub-affine complexity. 
With the help of the nice result of J. Cassaigne \cite{Ca} showing that a sequence $x$ has sub-affine complexity if and only if $(p_x(n+1) - p_x (n))_n$ is bounded, Ferenczi proved in \cite{Fe} that minimal subshifts with sub-affine complexity (i.e., $(p_X (n)/n)_n$ is bounded, where $p_X (n)$ is the number of words of length $n$ appearing in sequences of $X$) are $S$-adic (i.e.,  obtained by an infinite product of finitely many morphisms).
And, in the case it is ultimately bounded by $2$, he showed that less than $3^{27}$ morphisms are needed. 
The technicality of the proof did not allow to imagine and define the "strong notion of $S$-adicity" that is looked for.
In \cite{Le1} the second author presented a more detailed proof of Ferenczi's result that leaded to some improvements (see also \cite{LR}). Following Ferenczi's approach he showed that there exists a set $S$ consisting of 5 morphisms (to be compared to $3^{27}$) such that when $(p_X (n+1)-p_X (n))_n$ is bounded by $2$ then $X$ is described by an infinite product of morphisms belonging to $S$.
This provides a "strong notion of $S$-adicity" which is equivalent to this property.
When $(X,T)$ is such that $(p_X (n+1)-p_X (n))_n$ is bounded by $2$ for all large enough $n$, then another morphism is needed.  
We present this result in this note (Theorem \ref{theo:principal}) with a statement that applies to Bratteli-Vershik representations of such subshifts.
We recall that Bratteli-Vershik representations of minimal $\mathbb{Z}$-actions on Cantor sets are powerful tools that were used to characterize the topological orbit equivalence \cite{GPS} and that are very useful to solve problems where recurrence properties are involved. 
Observe it is usually not so easy to find such a representation and difficult to find the "canonical" representation. 
Nevertheless it has been done for various classical family of dynamical systems such that substitutive subshifts, Toeplitz subshifts, interval exchange transformations, sturmian subshifts, odometers and linearly recurrent subshifts.
In this note we present such representations for subshifts such that $(p_X (n+1)-p_X (n))_n$ is bounded by $2$ for all large enough $n$.

In Section 2 we recall the definition of Bratteli diagram and Bratteli-Vershik representation. 
Then we present some combinatorics on words and morphisms that notably help to obtain such representations of subshifts (Corollary \ref{coro:main}).
In Section 3, we present one of the main result of Leroy in \cite{Le2} and, as an application of this result, we obtain BV-representations of $S$-adic subshifts $(X,T)$ such that $(p_X (n+1)-p_X (n))_n$ is bounded by $2$ for all large enough $n$.

\section{Bratteli-Vershik representations of $S$-adic subshifts}
\label{sec:BV-Sadic}

\subsection{Bratteli-Vershik representations and subshifts}

Consider a minimal Cantor dynamical system $(X , T)$ (MCDS, also called minimal $\mathbb{Z}$-action on a Cantor set), i.e., a homeomorphism $T$ on a compact metric zero-dimensional space $X$ with no isolated points, such that the orbit
$\{ T^n x ; n \in \mathbb{Z} \}$ of every point $x\in X$ is dense in $X$. 
We assume familiarity of the
reader with the Bratteli-Vershik diagram representations of such systems, yet we
recall it briefly in order to establish the notation. 
For more details see \cite{HPS,Du}.
The vertices of the Bratteli-Vershik diagram $\mathcal{B} = (V,E)$ are organized into countably many
finite subsets of vertices $V_i $, $i\geq 0$ (where $V_0$ is a singleton $\{ v_0 \}$) and subsets of edges $E_i$, $i\geq 1$. 
Thus, $V=\cup_{i\geq 0} V_i$ and $E=\cup_{i\geq 1} E_i$. 
Every edge $e\in E_{i+1}$ connects
a vertex $s = s(e) \in  V_{i+1}$ for some $i\geq 0$ with some vertex $t = t(e) \in  V_i$.
At least one edge goes upward and at least one goes downward from each vertex in $V_{i+1}$.
Multiple arrows connecting the same vertices are admitted. 
We assume that the diagram is {\em simple}, i.e., there is a subsequence $(i_k)_{k\geq 0}$ such
that from every vertex in $V_{i_{k+1}}$
there is an upward path (going upward at each
level) to every vertex in $V_{i_k}$.
For each vertex $v\in V_{i}$ (except $v_0$) the set of all edges $\{ e_1, \dots , e_l \}$
going upward (to $V_{i-1}$) from $v$ is ordered linearly: let say $e_1< \dots < e_l$.
It will be convenient to consider the morphisms $\sigma_i^{\mathcal{B}} : V_i^*\to V_{i-1}^*$ defined by $\sigma_i^{\mathcal{B}} (v) = t(e_1) \cdots t(e_l)$.
We will say it is the morphism we read at level $i$ on $\mathcal{B}$.
In the sequel we will always suppose that $\mathcal{B}$ is such that $\sigma_1^{\mathcal{B}} (v)=v_0$ for all $v\in V_1$. 
This ordering induces a lexicographical order on all
upward paths from $v$ to $v_0$, and a partial order on all infinite upward paths arriving
to $v_0$. 
We denote by $X_\mathcal{B}$ the set of all such infinite paths. 
We assume that this
partial order has a unique minimal element $x_m$ (i.e., such that all its edges are
minimal for the local order) and a unique maximal one $x_M$ (i.e., such that all its edges are
minimal for the local order). 
This defines a map $V_\mathcal{B}$ sending every element $x$ to its successor in
the partial order and sending $x_M$ to $x_m$.

\begin{theorem}[\cite{HPS}]
Let $(X,T)$ be a minimal Cantor dynamical system. Then, there exists a simple Bratteli diagram $\mathcal{B}$ such that $(X,T)$ is conjugate to $(X_\mathcal{B} , V_\mathcal{B})$: there exists a bijective continuous map $\phi : X \to X_\mathcal{B}$ such that $\phi \circ T = V_\mathcal{B} \circ \phi$. 
\end{theorem}

This representation theorem has been used in \cite{GPS} to characterize the orbit equivalence relation of MCDS. 
We say that $(X_\mathcal{B} , V_\mathcal{B})$ is a BV-representation of $(X,T)$.

In the sequel $(X,T)$ will exclusively be a {\it subshift}: $X$ is a $T$-invariant closed subset of $A^{\mathbb{Z}}$ (endowed with the product topology) where $A$ is a finite alphabet and $T$ is the shift map ($T ((x_n)_n) = (x_{n+1})_n$). 
Note that the shift map will always be denoted by $T$: we will specify neither $X$ nor $A$.
We say that $(X,T)$ is generated by $x\in A^\mathbb{Z}$ when $X$ is the set of sequences $y$ such that for any $i$ and $j$ the word $y_i y_{i+1} y_{i+j}$ occurs in $x$. We set $[u.v]_X = \{ x \in X | x_{-|u|}\cdots x_{|v|-1} = uv \}$ and we call such sets {\em cylinder sets}. They are clopen sets and form a base of the topology of $X$. When $u$ is the empty word we write $[v]_X$.

\subsection{Proper morphisms, substitutions and BV-representations}

In the sequel $A$, $B$, $A_n$, $\dots$ are finite alphabets and $A^*$ denotes the free monoid generated by $A$.
A morphism $\sigma : A^* \to B^*$ is
{\em left proper} (resp. {\em right proper}) if there exists a letter $l\in B$ (resp. $r\in B$) such
that for all $a\in A$, $\sigma (a) = lu(a)$ (resp. $\sigma (a) = u(a) r$) for some $u (a)\in B^*$.
It is {\em proper} when it is both left and right proper.
Let $\sigma$ be left proper.
The morphism $\tau : A^* \to B^*$ defined by $\tau (a) = u(a)l$ is the {\em left conjugate} of $\sigma$ and it is right proper.
In the same way we define the {\em right conjugate} of $\sigma$ (it is left proper).

\begin{lemma}
\label{lemma:trick}
Let $\sigma : A^* \to B^*$ be a left proper (resp. right proper) morphism with first letter $l$ (resp. last letter $r$) and $\tau$ be its left (resp. right) conjugate then for all $a\in A$ and $n$

$$
\sigma^n (a)l = l\tau^n (a) \ \  (\hbox{resp. } r\sigma^n (a) =  \tau^n (a)r ).
$$
\end{lemma}

A substitution is an endomorphism $\sigma : A^* \to A^*$ such that there exists a letter $a$ with $\sigma (a) = au$ where $u$ is not the empty word. 
The subshift $(X_\sigma , T)$ it generates consists of $x\in A^\mathbb{Z}$ such that all words $x_ix_{i+1} \cdots x_j$ of $x$ have an occurrence in some $\sigma^n (a)$. 
We refer to \cite{Qu} for more details. 

To obtain a BV-representation of $(X_\sigma , T)$, an idea (first developed in \cite{Ve})
is to consider the Bratteli diagram $\mathcal{B}$ where the morphism we read at each level on $\mathcal{B}$ is $\sigma$.
When $(X_\sigma , T)$ is minimal this provides a measure-theoretical representation which is not necessarily topological. 
The problem with this construction is that minimal paths correspond to right fixed points and maximal paths to left fixed points of $\sigma$. 
Thus, to obtain such a topological representation, we need to have a unique fixed point in $A^\mathbb{Z}$. 

In \cite{DHS} is shown that to have a stationary BV-representation of $(X_\sigma , T)$ (i.e., where the substitutions read on the $E_i$ are equal up to some bijective changes of the alphabets), 
it suffices to find a proper substitution $\zeta$ such that $(X_\zeta , T)$ is conjugate to $(X_\sigma , T)$.
Then the stationary BV-representation is given by the Bratteli diagram where at each level $i\geq 2$ we read $\zeta$.
Moreover an algorithm is given to find such $\zeta$.
Consequently to Lemma \ref{lemma:trick} the following proposition claims it is enough to have $\zeta$ left or right proper.

\begin{proposition}
Let $\sigma$ be a left or right proper primitive substitution and $\tau $ be a conjugate. 
Then, 
\begin{enumerate}
\item
$\sigma\tau$ and $\tau \sigma$ are proper primitive substitutions;
\item
$(X_\sigma ,T)$, $(X_\tau , T)$, $(X_{\tau \sigma},T)$, $(X_{\sigma \tau },T)$, $(X_{\mathcal{B}_1},V_{\mathcal{B}_1})$ and $(X_{\mathcal{B}_2},V_{\mathcal{B}_2})$ are pairwise conjugate where $\mathcal{B}_1$ (resp. $\mathcal{B}_2$) is the stationary Bratteli diagram we read $\sigma \tau$ (resp. $\tau \sigma $) on from level $2$. 
\end{enumerate}
\end{proposition}

The algorithm used in \cite{DHS} to find BV-representations of substitution subshifts has a simplified version that can be used to find a left (or right) proper substitution. 
This leads to BV-representations with less edges and vertices. 
In the sequel we develop the idea of Lemma \ref{lemma:trick} to a more general framework.

\subsection{Combinatorics on words and BV-representations}

Let $S$ be a (possibly infinite) set of morphisms.
An {\em $S$-adic representation} of $(X,T)$ is a sequence $(\sigma_n , a_n )_{n\geq 2}$ where, for all $n$, $\sigma_n: A_{n}^* \rightarrow A_{n-1}^*$ belongs to $S$, $a_n$ belongs to $A_n$ and $X$ is the set of sequences $x\in A_1^\mathbb{Z}$ such that all words $x_ix_{i+1} \cdots x_j$ appear in some $\sigma_2 \sigma_3 \cdots \sigma_n (a_{n})$.
We start the representation with $n=2$ in order to fit with the Bratteli diagram notation: $\sigma_n$ will correspond to $E_n$ from $n=2$.
When such a sequence is fixed, we denote by $(X_n , T)$ the subshift generated by $(\sigma_k , a_k )_{k\geq n}$.
Observe $X_n$ is included in $A_{n-1}^\mathbb{Z}$.

\begin{proposition}
\label{prop:isom}
Let $(X,T)$ be the minimal $S$-adic subshift defined by $(\sigma_n : A_{n}^* \to A_{n-1}^* , a_n)_n$ where the $\sigma_n$ are proper.
Suppose that all morphisms $\sigma_n$ extend by concatenation to a one-to-one map from $X_{n}$ to $X_{n-1}$.
Then, $(X,T)$ is conjugate to $(X_\mathcal{B} , V_\mathcal{B})$ where $\mathcal{B}$ is the Bratteli diagram such that for all $n\geq 2$ 
the substitution read on $\mathcal{B}$ at level $n$ is $\sigma_{n}$.
\end{proposition}

{\em Proof.}
We can suppose that all the images of $\sigma_n$ starts with the letter $a$ and ends with $b$. 
Define, for all $n$, $\tau_n = \sigma_2 \cdots \sigma_n$ and $(X_n , T)$ the subshift generated by $(\sigma_l : A_{l}^* \to A_{l-1}^* , a_l)_{l\geq n}$.
Notice that $\tau_n (X_{n+1})$ is included in $X$ and consider

$$
{\mathcal P} (n) = 
\left\{
T^j \tau_n ([c]_{X_{n+1}})  \mid c \in A_{n} , 0\leq j< |\tau_n (c)|  
\right\} \ .
$$

From \cite{HPS}, to conclude, it suffices to prove that $({\mathcal P} (n))_n $ is a nested sequence of partitions generating the topology of $X$ satisfying $\# \cap_n\cup_{c \in A_{n}}\tau_n ([c]_{X_{n+1}}) = 1$.
We prove it is a partition. 
Let $x\in X$.
The maps $\sigma_n : X_{n} \to X_{n-1}$ being one-to-one, there exist a unique couple $(y,j)$, $y=(y_l)_{l\in \mathbb{Z}}\in X_{n}$ and $j\in \{ 0 ,1,\dots, |\tau_n (y_0)|-1 \}$, such that $x=T^j \tau_n (y)$. 
Then, $x$ belongs to $T^j \tau_n ([y_0]_{X_{n}})$ and ${\mathcal P} (n)$ is a partition. 
Now, we prove it is nested.
Let $\Omega = T^j \tau_{n+1} ([c]_{X_{n+2}})$ be an atom of $\mathcal{P} (n+1)$. 
Let $\sigma_{n+1} (c) = c_1 \cdots c_l$ and $i$ such that $|\tau_n (c_1 \cdots c_i )| \leq j < |\tau_n (c_1 \cdots c_{i+1} )|$. 
Then, $\Omega \subset T^{j- |\tau_n (c_1 \cdots c_i )|} \tau_{n} ([c_{i+1}]_{X_{n+1}})$
and $(\mathcal{P} (n))_n$ is nested.
We let as an exercise to prove it generates the topology of $X$.

As $\tau_{n+1} ([c]_{X_{n+2}}) \subset \tau_{n} ([b.a]_{X_{n+1}})$, from the assumptions, we deduce 
$\cap_n\cup_{c \in A_{n}}\tau_n ([c]_{X_{n+1}}) = 1$. 
\hfill $\Box$

Observe that in Proposition \ref{prop:isom}, it is important to suppose that $\sigma_n$ is one-to-one. 
For examples, if the $\sigma_n$ were equal to $\sigma$, where $\sigma (a)=ab$ and $\sigma (b)=ab$, then $(X,T)$ would be a two-point periodic subshift and $(X_\mathcal{B} , V_\mathcal{B})$ would be the $2$-odometer.
Using Lemma \ref{lemma:trick} we obtain the following corollary.

\begin{corollary}
\label{coro:read}
Let $(X,T)$ be the minimal $S$-adic subshift defined by $(\sigma_n , a_n)_{n\geq 2}$ where the $\sigma_n$ are left or right proper.
Suppose that for all $n$ the morphisms $\sigma_n$ extend by concatenation to one-to-one maps from $X_{n}$ to $X_{n-1} $.
Then  $(X,T)$ is conjugate to $(X_\mathcal{B} , V_\mathcal{B})$ where $\mathcal{B}$ is the Bratteli diagram where for all $n\geq 2$:
\begin{enumerate}
\item
the substitution read on $E_{2n}$ is left proper and equal to $\sigma_{2n}$ or its conjugate;
\item
the substitution read on $E_{2n+1}$ is right proper and equal to $\sigma_{2n+1}$ or its conjugate.
\end{enumerate}

\end{corollary}

More can be said about these subshifts but we need the following theorem proved in \cite{DM}.
We say that a MCDS $(Y,S)$ has {\em topological rank} $k$ if $k$ is the smallest integer such that $(Y,S)$ has a BV-representation $(X_\mathcal{B} , V_\mathcal{B})$ with $(|V (n)|)_n$ bounded by $k$. 
When such a $k$ does not exist we say that it has infinite topological rank.
Examples of such systems are the odometers ($k=1$), all primitive substitutive subshifts, all sturmian subshifts ($k=2$), all linearly recurrent subshifts, some toeplitz subshifts (see \cite{Du} for a survey). 

\begin{theorem}[\cite{DM}]
Let $(Y,S)$ be a MCDS with topological rank $k$.
Then $(Y , S)$ is expansive if, and only if, $k\geq 2$.
Otherwise it is equicontinuous.
\end{theorem}

In fact, we can be more precise.
Let $(X_\mathcal{B} , V_\mathcal{B} )$ be a Bratteli-Vershik representation of $(Y,S)$ whose set of vertices of $\mathcal{B}$ is bounded by $k$.
Now let $(X_{n} , T)$ be the subshift generated by $(\sigma_k^\mathcal{B} , a_k)_{k\geq n}$, where $a_k$ belongs to $V_{k}$.
Then, if there exists $n$ such that $(X_{n} , T)$ is not periodic then $( Y , S )$ is a subshift, otherwise it is equicontinuous. 
From this theorem and Corollary \ref{coro:read} we deduce the following corollary that enables us to give a BV-representation of any subshift once we have a "nice" $S$-adic representation.

\begin{corollary}
\label{coro:main}
Let $(X,T)$ be a MCDS with and $(\sigma_n : A_{n}^*\to A_{n-1}^* , a_n)_n$ be such that $(|A_n|)_n$ is bounded by $K$. 
The following are equivalent.

\begin{enumerate} 
\item
\label{cond:un}
$(X,T)$ is the non-periodic subshift defined by the sequence $(\sigma_n , a_n)_n$ satisfying: 
For all $n$ the map $\sigma_n : X_{n} \to X_{n-1} $ is well-defined, one-to-one and
left or right proper.
\item
\label{cond:deux}
$(X,T)$ is conjugate to $(X_\mathcal{B} , V_\mathcal{B})$ where $\mathcal{B}$ is the Bratteli diagram verifying:
\begin{enumerate}
\item
for all $n$, the substitution read on $E_{2n}$ is left proper and equal to $\sigma_{2n}$ or its conjugate;
\item
for all $n$, the substitution read on $E_{2n+1}$ is right proper and equal to $\sigma_{2n+1}$ or its conjugate;
\item
\label{cond:deuxc}
$(X_2 ,T)$, defined by $(\sigma_k , a_k)_{k\geq 2}$, is not periodic.
\end{enumerate}
\end{enumerate}
Moreover, in these situations the topological rank of $(X,T)$ is bounded by $K$.
\end{corollary}

In the next section we apply this corollary to the subshifts having a first difference of word complexity less or equal to 2 in order to obtain their Bratteli-Vershik representations.

\section{Applications to sub-affine complexity}

The following result was proven in \cite{Le1}. 
We adapt it to our context. 
The original statement gives a complete answer to the $S$-adic conjecture (see \cite{Fe}) in a restricted context.

\begin{theorem}
\label{theo:principal}
Let $(X,T)$ be an non-periodic minimal subshift such that $(p_X (n+1)-p_X (n))_n$ is bounded by $2$ for all large enough $n$, then
$(X,T)$ is an $S$-adic subshift defined by $(\sigma_n, a_n)_{n\geq 2}$ where for all $n \geq 3$, $\sigma_n$ belongs to $S = \{D,G,E_{ab},E_{bc},M \}$ with
$$
\begin{array}{llllllllllllll}
D:   & a\mapsto ab, & & G:   & a \mapsto ba, & & E_{ab}:& a \mapsto b, & & E_{bc}: & a \mapsto a, & & M: & a \mapsto a \\
    & b \mapsto b  & &     & b \mapsto b   & &       & b \mapsto a  & &        & b \mapsto c  & &   & b \mapsto b \\
    & c \mapsto c  & &     & c \mapsto c   & &       & c \mapsto c  & &        & c \mapsto b  & &   & c \mapsto b. 
\end{array}
$$
Moreover, $\sigma_n : X_{n} \to X_{n-1}$ is one-to-one and there exists an increasing sequence $(n_i)_i$ such that for all $i$, $\sigma_{n_i} \sigma_{n_i + 1} \cdots \sigma_{n_{i+1}}$ is proper and all letters in $\{a,b,c\}$ occur in all images $\sigma_{n_i} \sigma_{n_i + 1} \cdots \sigma_{n_{i+1}}(y)$, $y \in \{a,b,c\}$. 
Furthermore, $\sigma_2$ can be algorithmically determined.
\end{theorem}

Theorem~\ref{theo:principal} is weaker than the one in~\cite{Le1}. 
Indeed, the complete result states that there is a labelled directed graph $\mathcal{G}$ and some condition (C), both computable, for which a subshift $(X,T)$ is such that $(p_X (n+1)-p_X (n))_n$ is bounded by $2$, for all large enough $n$, if and only if it has a $S$-adic representation $(\sigma_n)_n$ which is an infinite path in $\mathcal{G}$ satisfying the condition (C). 

The Bratteli-Vershik version of this result is the following. 

\begin{corollary}
With the assumptions and notations of Theorem \ref{theo:principal}, the topological rank of $(X,T)$ is at most $3$.  
More precisely, $(X,T)$ is conjugate to $(X_\mathcal{B} , V_\mathcal{B})$ where $\mathcal{B}$ is the Bratteli diagram such that, for all $n\geq 3$, the substitution read at level $n$ is $\sigma_n$. 
\end{corollary}


\begin{thebibliography}{00}
\bibitem{AR}
P. Arnoux, G. Rauzy,
Repr\'esentation g\'eom\'etrique de 
suites de complexit\'e $2n+1$,
Bull. Soc. Math. France 119 (1991), 199-215.

\bibitem{Ca}
J. Cassaigne,
Special factors of sequences with linear subword complexity,
Developments in language theory II (Magdeburd, 1995), 25-34, World Sci. Publ., RiverEge, NJ, 1996.

\bibitem{DM}
T. Downarowicz, A. Maass,
Finite rank Bratteli-Vershik diagrams are expansive,
Ergod. Th. \& Dynam. Sys. 28 (2008), 739-747.

\bibitem{Du}
F. Durand,
Combinatorics on Bratteli diagrams and dynamical systems, 
Combinatorics, Automata and Number Theory, Series Encyclopedia of Mathematics and its applications 135, Cambridge University Press 2010, 338-386.

\bibitem{DHS} 
F. Durand, B. Host, C. Skau, 
Substitutive dynamical systems, Bratteli diagrams and dimension 
groups, Ergod. Th. \& Dynam. Sys. 19 (1999), 953-993. 

\bibitem{Fe}
S. Ferenczi,
Rank and symbolic complexity,
Ergod. Th. \& Dynam. Sys. 16 (1996), 663-682.

\bibitem{GPS}
T. Giordano, I. Putnam, C. Skau,
Topological orbit equivalence and $C^{*}$-crossed products,
Internat. J. Math. 469 (1995),
51-111.

\bibitem{HPS} 
R. H. Herman, I. Putnam, C. F. Skau, 
Ordered Bratteli diagrams, dimension groups and topological 
dynamics, Internat. J. of Math. 3 (1992), 827-864.

\bibitem{Le1}
J. Leroy, 
Some improvements of the $S$-adic conjecture, 
Adv. in Applied Math. 48 (2012), 79-98.

\bibitem{Le2}
J. Leroy, Contribution à la r\'esolution de la conjecture S-adique,
PhD thesis, Univ. Picardie Jules Verne (2011).


\bibitem{LR}
J. Leroy, G. Richomme, A combinatorial proof of $S$-adicity for sequences with sub-affine complexity, (preprint).


\bibitem{MH1}
M. Morse, G. A. Hedlund,
Symbolic dynamics,
Amer. J. Math. 60 (1938), 815-866.

\bibitem{MH2}
M. Morse, G. A. Hedlund,
Symbolic Dynamics {II}.  {Sturmian} trajectories,
Amer. J. Math. 62 (1940), 1-42.

\bibitem{Qu} 
M. Queff\'elec, 
Substitution Dynamical Systems-Spectral Analysis, 
Lecture Notes in Mathematics, 1294, Springer-Verlag, Berlin, 1987.

\bibitem{Ve}
A. M. Vershik,
A theorem on the Markov periodical approximation in ergodic theory,
J. Sov. Math. 28 (1985), 667-674.

\end{thebibliography}
\end{document}